\documentclass[11pt]{article}

\usepackage{amsmath,amssymb,graphicx,dutchcal}

\def\mr{$^{\text{r}}$ }

\title{Hermite's letters to Francisco Gomes Teixeira}

\author{Pedro J.\ Freitas\\ Centro Interuniversit\'ario de Hist\'oria das Ci\^encias e da Tecnologia e \\Departamento de Hist\'oria e Filosofia das Ci\^encias,\\ Faculdade de Ci\^encias, Universidade de Lisboa,\\ 1749-016 Lisboa, Portugal\\
\\pjfreitas@fc.ul.pt}

\begin{document}

\maketitle 
\thispagestyle{empty} 

\setcounter{page}{1}

\centerline{\Large Hermite's letters to Francisco Gomes Teixeira}

\begin{abstract}
It is well known that Charles Hermite kept an intense correspondence with many of the word's leading mathematicians of his time. This paper focuses on Hermite's letters to Francisco Gomes Teixeira, a Portuguese mathematician, who exchanged letters with Hermite for twenty years. 
\end{abstract}

\section{Introduction}

Francisco Gomes Teixeira (1851-1933) was one of the most important Portuguese mathematicians of the late 19th century and early 20th. Having graduated from the University of Coimbra in 1875, he got his doctorate later in the same year, and, in the next year, started teaching at this university. A few years later, in 1884, he moved to the Polytechnic Academy of Porto, which would later become the University of Porto, in 1911. At this time, he became Rector of the newly established university. 

In addition to his scientific works, mostly in analysis, published in both Portuguese and foreign journals, he innovated in the field of teaching, by writing college manuals with ambitious content, which elevated the level of rigour of his time, in Portugal.\footnote{The manuals were {\it Curso de An\'alise infinitesimal, c\'alculo diferencial}, Porto, Typ.\ Ocidental, 1887 and {\it C\'alculo integral}, Porto: Typ.\ Occidental, 1889, edited and extended in later editions.}

One may say that Gomes Teixeira's most relevant contribution for the evolution of mathematics in Portugal was the broadening of its scope, from the national to the international level. He received two prizes from the Royal Academy of Sciences of Madrid, in 1895 and 1897, which are representative of this interest in establishing contacts at an international level. The second of these prizes was attributed to a treatise on special curves, {\em Tratado de las curvas especiales notables, tanto planas como alabeadas}, which was later expanded, translated to French, and republished in 1917, after receiving another prize from the French Academy of Sciences. For more on his life and work, see for instance \cite{Vi}, \cite{Al}.\bigskip

Along with the activity we have mentioned, there were several aspects which are central to this internationalization of Portuguese mathematics achieved by Gomes Teixeira. One of these was the foundation, in 1877, of the scientific journal {\it Jornal de Sciencias Matem\'aticas e Astron\'omicas}, which we will simply refer to as JSMA. For more on this journal, see \cite{Sa} and \cite{Ka}, a thesis that studies importance of the journal in its time.\footnote{The tomes of this journal are now fully digitized and can be found at {\tt www.fc.up.pt/fa/index.php?p=nav\&f=html.fbib-Periodico-oa.}}

The Journal benefited greatly from Gomes Teixeira's intense correspondence with some of the most renowned mathematicians of his time. The letters that Gomes Teixeira received are kept in the Archive of the University of Coimbra, which includes more that two thousand letters, indexed and catalogued by Gomes Teixeira himself. As a contribution to the analysis of this estate, this paper focuses on his correspondence with Charles Hermite, one of the most represented authors: there are 19 letters from Hermite, some of them photographically reproduced in \cite{Al} and \cite{Ka}, which we present here, and comment for the first time.

\section{The correspondence}

Charles Hermite (1822--1901) was known for being a very prolific correspondent. Goldstein's paper \cite{Go} gives a general overview of Hermite's abundant correspondence (thousands of letters written). The paper also notes that since the beginning of the 20th century, the letters Hermite sent to many mathematicians have been gradually edited---our present paper wishes to continue this trend. 

Most of the correspondence received by Hermite was lost in a fire, which unfortunately makes it impossible to confront the letters in the Coimbra archive with the ones sent by Gomes Teixeira. The letters we present here go from 1872 to 1896, the number of letters for each year is as follows.

\begin{center}
\begin{tabular}{|l|c|c|c|c|c|c|c|c|c|c|c|c|}
\hline 
Year & 72 & 75 & 79 & 81 & 84 & 85 & 86 & 88 & 90 & 91 & 92 & 96\\ \hline
No. of letters & 1 & 2 & 1 & 2 & 1 & 2 & 2 & 2 & 2 & 1 & 2 & 1\\ \hline
\end{tabular} 
\end{center}


There are nine other letters referring to the homage organized for Hermite on the occasion of his 70th birthday, in 1892---a subscription was made to present Hermite with a gold medal struck in his honor. Gomes Teixeira was invited by Mittag-Leffler to be a member of the honor committee for this event. The letters are from Mittag-Leffler, Darboux, M.\ Merino, Loriga, Galdeano, C. Ricart and other undisclosed senders, most of them asking how they participate in the homage.

\subsection{July 17, 1872}

This is letter 141, and is the earliest letter from Hermite that can be found in this archive. Gomes Teixeira was still a university student at this point, he was finishing the 3rd year and starting the 4th. The presence of this letter, with a clear reference to the JSMA he would establish five years later, reveals that, even as a 21 year-old student, Gomes Teixeira already envisaged a new Portuguese mathematical journal and did not hesitate to write to one of the most famous mathematicians of his time asking for collaboration. 

The paper contained in the letter is \cite{HJ1}, on the Frenet-Serret formulas for curves in 3-dimensional space. It was the first paper by a non-Portuguese author published in the JSMA. It can also be found in \cite[p. 508-511]{HOC}.

\begin{quote}
Monsieur, 

Vous m'avez demand\'e, en commen\c cant la publication de votre Journal de Sciences Math\'ematiques et astronomiques, de vous donner ma collaboration; je viens remplir l'engagement que j'ai pris envers vous en vous donnant la note ci-jointe, concernant un [point?] \'el\'ementaire du calcul diff\'erentiel. 

Veuillez agr\'eer Monsieur l'expression de mes sentiments les plus distingu\'ees, 

Ch.\ Hermite 

Paris 17 Julliet 1872
\end{quote}

\subsection{June 10, 1875}

This is letter 140 in the archive. In this letter, Hermite thanks Gomes Teixeira for sending him his inaugural dissertation (this is the name used for the Ph.D.\ thesis, see \cite[pp.\ 34ff.]{Al}), entitled {\em Integra\c c\~ao das equa\c c\~oes \`as derivadas parciais de 2a.\ ordem} (Integration of 2nd order partial differential equations). He takes the opportunity to send some publications, without reference to its content. 

\begin{quote}
Monsieur, 

Je viens de parcourir la th\`ese inaugurale que vous m'avez fait l'honneur de m'adresser et quoique le sujet important et difficile que vous avez trait\'e ne reste point dans le cercle habituel de mes \'etudes, j'ai acquis l'assurance que vous avez fait un travail tr\`es s\'erieux et approfondi. 

Veuillez accepter, Monsieur, les opuscules qui accompagnent cette lettre comme un t\'emoignage de ma sympathie et l'expression de ma consid\'eration distingu\'ee.

Ch.\ Hermite

Paris, 10 juin 1875
\end{quote} 

\subsection{December 3, 1875}

This is letter 765 in the archive, and contains an exercise given by Hermite to his students as a competition, on the subject of continued fractions, one of Teixeira's interests. The exercise has some interest outside the context of student competitions: it provides a new proof to a result by Lambert, quoted in the letter. Lambert is credited with the first proof that $\pi$ is irrational (a proof which also uses continued fractions) in \cite{Lb}. This paper also studies powers of $e$, so this is probably the paper Hermite refers to in the letter.\footnote{Hermite himself proved, in 1873, that $e$ is transcendental, the result appears in two letters, sent to Gordan and Borchardt.}

We note that the exercise in this letter, with its resolution, does not appear among Hermite's contributions to the {\it Nouvelles Annalles de Math\'e\-ma\-tiques}, nor is it published in Gomes Teixeira's JSMA. 

\begin{quote}
3 D\'ec 1875

Monsieur, \medskip

Je saisis Monsieur cette occasion pour vous donner communication d'une petite question qui a \'et\'e le sujet de la derni\`ere comp\'etition donn\'e a mes \'el\`eves et qui se rapporte pr\'ecis\'ement a la th\'eorie des fraction continues. J'ai propos\'e de d\'emontrer qu'en posant 
$$F(x) = 1+ \frac x1 + \frac{x^2}{1.2} + \ldots + \frac{x^n}{1.2.\ldots.n}$$
de sorte que 
$$e^x = F(x) + \frac{x^{n+1}}{1.2.\ldots.(n+1)} + \text{etc\ldots}$$
la d\'eriv\'e d'ordre $n$ de la expression 
$$\frac{e^x-F(x)}{x^{n+1}}$$
a la forme suivante:
$$\frac{e^x\varphi(x) -\psi(x)}{x^{2n+1}}$$
o\`u $\varphi(x)$ et $\psi(x)$ sont des polyn\^omes, \`a coefficients entiers, du degr\'e $n$. Vous voyez ainsi que $\frac{\psi(x)}{\varphi(x)}$ est l'une des r\'eduites de la fraction continue qui repr\'esente l'exponentielle $e^x$, et cette circonstance que les coefficients de $\varphi(x)$ et $\psi(x)$ sont des nombres entiers conduit a une d\'emonstration immediate de la proposition d\'ecouverte par Lambert, que toutes les puissances de $e$ sont incommensurables. Il suffit en effet d'obtenir que la d\'eriv\'e d'ordre $n$ de la s\'erie 
$$\frac1{1.2.\ldots . (x+1)} + \frac{x}{1.2.\ldots . (x+2)} + \cdots$$
\'etant:
\begin{eqnarray*}
\frac1{(n+1).(n+2)\ldots (2n+1)} & + & \frac1{(n+2).(n+3).\ldots . (2n+2)}\, \frac{x}1\\
& + & \frac1{(n+3)(n+4)\ldots (2n+3)}\, \frac{x^2}{1.2} \\
& + & \frac1{(n+4)(n+5)\ldots (2n+4)}\,  \frac{x^3}{1.2.3}  \\
& + & \cdots
\end{eqnarray*}
on a pour la quantit\'e:
$$e^x\varphi(x) - \psi(x)$$
un d\'eveloppement en s\'erie qu'on peut \'ecrire sous la forme suivante:
$$\frac{x^{2n+1}}{(n+1)(n+2)\ldots(2n+1)} \left[ 1+\frac{n+1}{2n+2} \, \frac{x}1 + \frac{(n+1)(n+2)}{(2n+2)(2n+3)} \, \frac{x^2}{1.2} +\cdots \right]$$
On voit ainsi que le nombre de $n$ croissant, cette expression peut devenir sans jamais s'annuler, plus petite que toute quantit\'e donn\'ee; car le facteur 
$$\frac{x^{2n+1}}{(n+1)(n+2)\ldots(2n+1)}$$
a pour limite, z\'ero, tandis que la s\'erie entre crochets est toujours finie, les termes \'etant respectivement moindres que ceux de $e^x = 1+\frac{x}1 + \cdots$

Veuillez Monsieur m'excuser du retard que mes occupations m'obligent de mettre \`a la lecture de votre travail et recevoir la nouvelle assurance de ma consid\'eration distingu\'ee.

Ch.\ Hermite 

Paris, 2 Rue da la Sorbonne, 3 D\'ecembre 1875.

\end{quote}

\subsection{November 19, 1879}

This is letter 142 in the archive, and contains the paper \cite{HJ2}, which can also be found at \cite[pp. 505--507]{HOC}. It is a result about an integral involving trigonometric functions, which Hermite proves in his {\it Cours d'Analyse}. The paper contains a simpler proof of the result. The letter also contains some results about interpolating polynomials.

\begin{quote}

Paris, 19 Novembre 1879

(2 Rue de Sorbonne)

Monsieur,

Ce m'est un grand regret que sur le petit nombre de 20 exemplaires de tirage \`a part de mon [article?] sur l'article concernant l'interpolation, il ne m'en reste plus un seul que je puisse vous offrir pour r\'epondre \`a votre d\'esir. Permettez moi de me d\'edommager un tant soit peu en vous envoyant avec cette lettre toutes celles de mes publications dont je puis disposer et aussi en vous donnant en quelques mots le principe de ma formule. La question \'etant de d\'eterminer un polyn\^ome $F(x)$ de degr\'e $n-1$ satisfaisant aux conditions:
\begin{eqnarray*}
&& F(a)=f(a), F(b)=f(b), \ldots F(l) = f(l) \\
&& F'(a)=f'(a), F'(b)=f'(b), \ldots F'(l) = f'(l) \\
&& \ \ \vdots \\
&& F^{\alpha-1}(a)=f^{\alpha-1}(a), F^{\beta-1} (b)=f^{\beta-1}(b), \ldots F^{\lambda-1}(l) = f^{\lambda-1}(l) 
\end{eqnarray*}
ou $f(x)$ est une fonction donn\'ee; en supposant $\alpha+\beta+\ldots+\lambda=n$, je consid\`ere une aire $S$, comprenant d'une part, $a$, $b$, \ldots, $l$, et de l'autre la quantit\'e variable $x$, j'admets qu'\`a son int\'erieur la fonction $f(x)$ soit uniforme et n'ait aucun p\^ole; cela \'etant, on a:
\begin{eqnarray*}
&& F(x)-f(x) =\\
&& \quad = \frac1{2i\pi} \int \frac{f(z).(x-a)^\alpha (x-b)^\beta \ldots (x-l)^\lambda}{(x-z).(z-a)^\alpha (z-b)^\beta \ldots (z-l)^\lambda} dz
\end{eqnarray*}
l'int\'egrale du second membre se rapportant au contour de $S$. L'utilit\'e principale de cette relation est moins de d\'eterminer le polyn\^ome cherch\'e $F(x)$ par le calcul de l'int\'egrale, et au moyen d'une somme de r\'esidus, que de montrer que la diff\'erence $F(x)-f(x)$, diminue sans limite, lorsque le nombre des quantit\'es, $a$, $b$, \ldots, $l$, ou bien les exposants $\alpha$, $\beta$, \ldots, $\lambda$ vont en augmentant. 

A ma bien l\'eg\`ere offrande d'une note pour votre journal, je [joints?] Monsieur l'expression de ma consid\'eration la plus distingu\'ee et de mes sentiments d\'evou\'es

Ch.\ Hermite

\end{quote}

\subsection{July 16, 1881}

This is letter 143 in the archive. In this letter, Hermite comments on a note Gomes Teixeira had sent him for presentation at the Acad\'emie des Sciences, stating that this result was already known. About this letter, Gomes Teixeira writes in the catalogue (see \cite{Vi}) that he later sent a new result which was indeed published in the {\it Comptes Rendues}, namely \cite{GT2}. The topic of this second note is the integration of a partial differential equation, the same topic of his thesis. Letter 764 (immediately below) acknowledges the reception of this second paper, and contains some more comments. 

It is remarkable that, after gracefully declining the publication of the first note, Hermite encourages Gomes Teixeira to send in a new one, acknowledging him as a representative of mathematics both for Portugal and Spain.

\begin{quote}

Paris, 16 Juillet 1881

Monsieur

Le travail que vous m'avez communiqu\'e rapportant a des questions d'analyse dont M\mr Darboux s'est occup\'e sp\'ecialement et avec le plus grand succ\`es, je l'ai pri\'e de vouloir bien en prendre connaissance et c'est son avis autant que le mien dont je viens vous faire part.

Nous croyons qu'il vous serait n\'ecessaire d'\'etudier plusieurs recher\-ches r\'ecentes et importantes dont vous ne paraissez pas avoir eu connaissance, notamment celle de Bour\footnote{Edmond Bour (1832--1866) studied both at the \'Ecole Polyt\'echnique and the \'Ecole des Mines. He published on differential equations, celestial mechanics and geometry of surfaces.} et aussi les m\'emoires 
d'[Imenetsky?] et de M\mr Pelet.\footnote{This may be a reference to Auguste Pellet (1848--1935), who worked mostly on Analysis, from 1867 until the 1910's.} 

Votre nom repr\'esentant les sciences math\'ematiques pour votre pays et aussi pour l'Espagne, nous croyons devoir vous conseiller de ne le produire pour la premi\`ere fois dans les Comptes-Rendues qu'a l'occasion d'un travail compl\`etement digne de votre situation. Votre talent Monsieur, et je suis heureux de vous le dire au nom de M\mr Darboux comme au mien, vous appelle a le remplir de la mani\`ere la plus honorable. Nous pensons mieux vous t\'emoigner notre sympathie en vous demandant de nouveaux efforts qu'en publiant un travail dont on pourrait dire qu'il n'a pas \'et\'e accueilli sans quelque indulgence, et c'est en esp\'erant que nous aurons bient\^ot le fruit de ces efforts que je vous offre Monsieur la nouvelle expression de toute mon estime et de mes sentiments bien d\'evou\'es, 

Ch.\ Hermite 

\end{quote}

\subsection{November 1, 1881}

This is letter 764 in the archive. 

\begin{quote}

Monsieur, 

Je m'empresse de vous accuser r\'eception de la note que vous m'avez adress\'ee et de vous informer que je remplirai vos intentions avec grand plaisir en la pr\'esentant a la prochaine S\'eance de l'Acad\'emie, afin qu'elle soit publi\'ee dans les Comptes-Rendues. J'esp\`ere que plutard\footnote{This word is not currently in use, but it does appear in other 19th century sources, e.g., in a letter from Flaubert to George Sand, dated from July 27, 1871---see \cite{La}.} vous donnerez dans votre journal des applications qui en feront ressortir l'utilit\'e et l'importance.

Veuillez agr\'eer Monsieur la nouvelle assurance de mes sentiments bien d\'evou\'es, 

Ch. Hermite

Paris 1 Novembre 1881. 

\end{quote}

\subsection{December 14, [1884?]}

This is letter 769 in the archive. In this letter there is no reference to any paper. Apparently, Gomes Teixeira sent a theorem to Hermite (possibly in a previous letter), without proof, and Hermite encourages him to find one, taking the opportunity to make interesting remarks about the need for the presentation of mathematical results to be as simple and clear as possible, a constant concern in these letters. 

\begin{quote}
Paris 14 D\'ecembre [1884?]

Monsieur, 

L'\'enonc\'e du th\'eor\`eme que vous m'avez fait l'honneur de m'adresser me fait voir que vous \^etes engag\'e dans une voie excellente et je ne puis que vous engager vivement \`a faire tous les efforts n\'ecessaires pour arriver a une d\'emonstration simple et facile que d'elle-m\^eme devienne classique. 
Vous le savez Monsieur, et je n'ai pas a vous l'apprendre qu'apr\`es le travail de l'invention, il en est un autre dont nos ma\^itres en Analyse ont donn\'e l'exemple et le mod\`ele, de sorte que les \oe uvres 
de Gauss et de Jacobi joignent \`a l'importance de l'\'eclat des d\'ecouvertes, le m\'erite d'une forme parfaite. En vous exprimant l'espoir qu'ayant \'et\'e assez heureux pour obtenir un beau r\'esultat, vous [?]\footnote{Illegible.} aussi l'exposer sous la forme la meilleure, je saisis Monsieur cette occasion pour vous renouveler l'expression de tous mes sentiments de sympathie et d'estime.

Ch.\ Hermite.
\end{quote}

\subsection{May 31, 1885}

This is letter 766 in the archive and contains the paper \cite{HJ3}, which also appears in \cite[pp.\ 169--171]{HOC}. It concerns a relation between Legendre polynomials and continued fractions of functions (we include the beginning of the paper in the transcription below). 

\begin{quote}

Paris, 31 Mai 1885

Monsieur,

Je viens vous remercier les explications que vous avez eu la bont\'e de me donner et qui m'ont permis de mieux saisir votre analogie, \`a laquelle je n'ai plein d'objections \`a faire. Les applications de la formule \`a laquelle vous \^etes parvenu, \`a des cas particuliers simples me semblait utiles pour en faire comprendre le caract\`ere; mais en ce moment je ne puis m\^eme songer a vous donner d'indications pr\'ecises comme il serait n\'ecessaire, les devoirs dont je suis oblig\'e ne me laissant aucune libert\'e d'esprit.

Je voudrais cependant Monsieur vous offrir un t\'emoignage de mes sentiments de sympathie et d'estime en vous donnant pour \^etre publi\'e dans votre journal, si vous le voulez bien, une petite remarque qui a fait l'objet d'une de mes derni\`eres le\c cons \`a la Sorbonne. 

Vous connaissez cette belle proposition de M. Tchebyshew  que le polyn\^ome $X_n$ de Legendre est le d\'enominateur de la r\'eduite d'ordre $n$ du d\'eveloppent  en fraction continue de la quantit\'e:
$$\frac12 \log \frac{x+1}{x-1} = \frac1x + \frac1{3x^3} + \frac1{5x^5} + \cdots$$

On peut y parvenir comme vous allez voir, au moyen du d\'eveloppement en s\'erie qui a \'et\'e le point de d\'epart de Legendre et a donn\'e la premi\`ere des formules des polyn\^omes $X_n$. [...]

Veuillez Monsieur recevoir la nouvelle assurance de toute ma sympathie et de mes sentiments bien d\'evou\'es

Ch.\ Hermite
\end{quote}

\subsection{October 20, 1885}

This is letter 767 in the archive. The paper mentioned in this letter is \cite{GT3}, which contains a proof of Eisenstein's theorem on power series\footnote{See the wikipedia page {\tt http://en.wikipedia.org/wiki/Eisenstein's$\_$theorem}\ for the statement of the theorem.} (it proves a more general statement). The next letter, from January 26, 1886 (number 770), refers to the same paper. 

\begin{quote}

Monsieur, 

J'esp\`ere n'avoir pas \'et\'e contre vos intentions en demandant \`a mon confr\`ere M\mr Darboux de publier dans les Annales de l'\'Ecole Normale Sup\'erieure le travail que vous avez fait l'honneur de me communiquer. La propri\'et\'e que vous avez ajout\'ee aux r\'esultats importants d\'ecouverts par Eisenstein sur les s\'eries a coefficients rationnels qui satisfont a une \'equation alg\'ebrique a coefficients entiers nous a tous deux beaucoup int\'eress\'es, et je serais heureux Monsieur de joindre des f\'elicitations \`a l'assurance de ma haute estime et \`a celle de mes sentiments bien d\'evou\'es.

Ch.\ Hermite 

20 Octobre 1885
\end{quote}

\subsection{January 26, 1886}

This is letter 770 in the archive. It refers to the same paper as the previous one, with a special note to the elegance and simplicity of the proof provided.

\begin{quote}

Monsieur

J'ai donn\'e communication a M\mr Darboux de la d\'emonstration du th\'eor\`eme d'Eisenstein que vous m'avez fait l'honneur de m'\'ecrire, et nous pensons remplir vos intentions en la publiant dans les Annales de l'\'Ecole Normale, \`a la suite du premier article que vous avez donn\'e a ce recueil.

Sinc\`eres f\'elicitations Monsieur de la simplicit\'e et de l'\'el\'egance de votre d\'emonstration, je vous prie d'agr\'eer l'assurance de toute mon estime et de mes sentiments d\'evou\'es, 

Ch.\ Hermite 

Paris 26 Janvier 1886

\end{quote}

\subsection{February 12, 1886}

This is letter 768 in the archive. The article this letter refers to is \cite{GT4}, which is a development of the result in paper \cite{GT3}, published the year before, also in the {\em Annales de l'\'Ecole Normale Sup\'erieure}.

\begin{quote}
Monsieur,

J'ai le plaisir de vous informer que la derni\`ere communication que vous avez fait l'honneur de m'adresser sera publi\'e, comme les pr\'ec\'edentes, dans les Annales de l'\'Ecole Normale.

En vous renouvelant Monsieur mes sinc\`eres f\'elicitations pour ce nouveau r\'esultat de votre travail, je vous prie d'agr\'eer l'assurance de mes sentiments les plus distingu\'es.

Ch.\ Hermite 

Paris 12 F\'evrier 1886
\end{quote}

\subsection{November 16, 1888}

This is letter 772 in the archive. It mentions two letters sent by Gomes Teixeira. In the present letter, Hermite provides a new proof of a formula sent in the second of Gomes Teixeira's letters. He also offers publication of a paper, probably \cite{GT51} (it may have been sent with the first letter). 

There are two figures illustrating the text, we include the original ones along with a reproduction. 

\begin{quote}
10 Novembre 1888

En lisant avec attention la seconde des deux lettres que vous m'avez fait l'honneur de m'adresser et dans laquelle vous parvenez a cette relation:
\begin{eqnarray*}
&& \int_a^b \varphi^2(x)dx.\!\int_a^b \psi^2(x)dx - \left[  \int_a^b \varphi(x)\psi(x)dx \right]^2 \\
&&\qquad = 
\int_a^b dx \int_x^b [\varphi(x)\psi(y) -\varphi(y)\psi(x)]^2 dy
\end{eqnarray*}
ou bien, dans une forme un peu plus g\'en\'erale, 
\begin{eqnarray*}
&& \int_a^b \varphi(x)\varphi_1(x)dx.\!\int_a^b \psi(x)\psi_1(x)dx -  \int_a^b \varphi(x)\psi_1(x)dx.\int_a^b \varphi_1(x)\psi(x)dx  \\
&&\qquad = 
\int_a^b dx \int_x^b [\varphi(x)\psi(y) - \varphi(y) \psi(x)] 
[ \varphi_1(x)\psi_1(y) - \varphi_1(y) \psi_1(x)]  dy,
\end{eqnarray*}
j'ai reconnu que cette relation est identique comme vous allez le voir.

Consid\'erons en effet l'int\'egrale double:
$$J=\int_a^x dx \int_x^b F(x,y)dy.$$
Il repr\'esente le volume d'une section du prisme limit\'e par la surface $z=F(x,y)$ et compris entre trois plans d\'etermin\'es comme il suit.

\begin{center}
\begin{picture}(250,100)
\put(130,0){\includegraphics[scale=0.15]{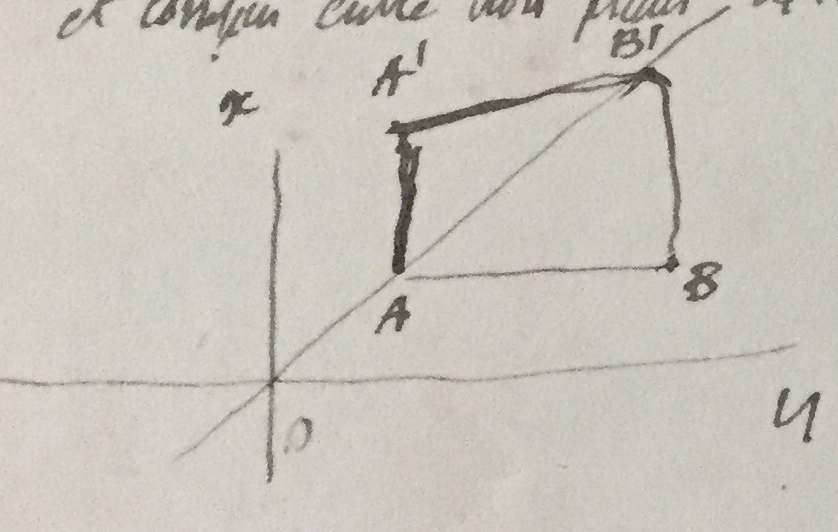}}
{\fontsize{10pt}{10pt}
\put(0,10){\line(1,0){100}}
\put(10,0){\line(0,1){100}}
\put(0,0){\line(1,1){100}}
\put(25,25){\line(1,0){60}}
\put(25,25){\line(0,1){60}}
\put(85,85){\line(-1,0){60}}
\put(85,85){\line(0,-1){60}}
\put(22,14){$A$}\put(82,14){$B$}
\put(22,88){$A'$}\put(78,88){$B'$}
\put(13,2){$o$}\put(2,90){$x$}\put(90,3){$y$}
}
\end{picture}
\end{center} 

Tracez sur le plan des $xoy$ la bissectrice $y=x$ de l'angle $xoy$ puis le rectangle $ABB'A'$ dont elle serait la diagonale, et supposons que les abcisses des sommets $A$ et $B$ du rectangle soient $a$ et $b$. Les plans dont je parle, perpendiculaires au plan des $xy$ auront pour 
traces les droites, $AA'$, $A'B'$, $AB'$ et la base du prisme consid\'er\'e sur le triangle rectangle $AA'B$. 

\begin{center}
\begin{picture}(250,100)
\put(130,0){\includegraphics[scale=0.15]{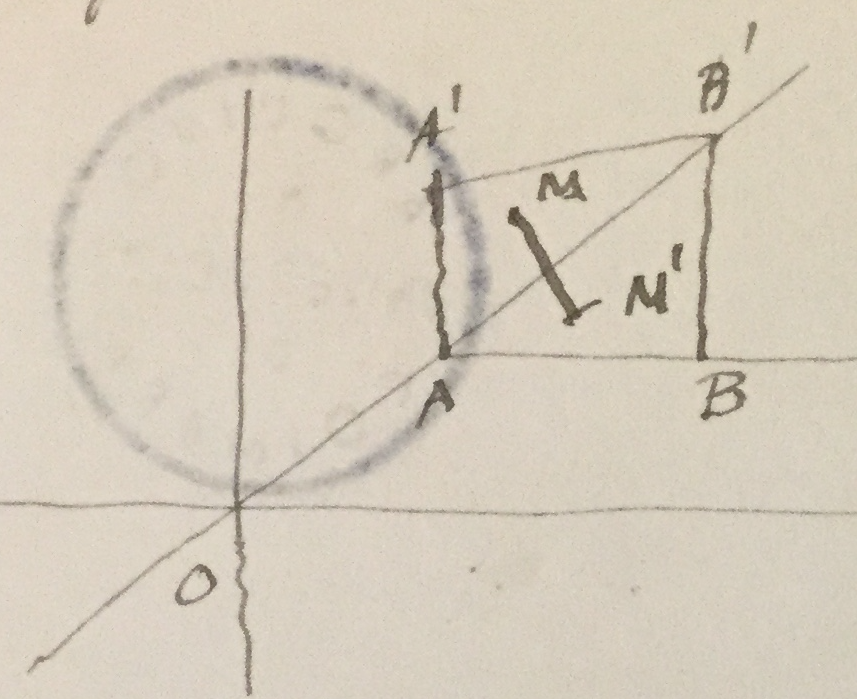}}
{\fontsize{10pt}{10pt}
\put(0,10){\line(1,0){100}}
\put(10,0){\line(0,1){100}}
\put(0,0){\line(1,1){100}}
\put(25,25){\line(1,0){60}}
\put(25,25){\line(0,1){60}}
\put(85,85){\line(-1,0){60}}
\put(85,85){\line(0,-1){60}}
\put(40,60){\line(1,-1){20}}
\put(40,60){\circle*{2}}\put(60,40){\circle*{2}}
\put(35,64){$M$}\put(63,35){$M'$}
\put(22,14){$A$}\put(82,14){$B$}
\put(22,88){$A'$}\put(78,88){$B'$}
\put(13,2){$o$}
}
\end{picture}
\end{center} 

Soit $M$ un point de l'int\'erieur de ce triangle, ayant pour coordonn\'ees $x$ et $y$, son sym\'etrique $M'$ par rapport a $AB'$ aura \'evidemment pour coordonn\'ees $y$ et $x$, par cons\'equent les valeurs de l'ordonn\'ee $z=F(x,y)$, seront les m\^emes en ces deux points lorsque la fonction $F(x,y)$ est sym\'etrique en $x$ et $y$. C'est le cas de l'expression qui entre dans votre int\'egrale double; on peut par cons\'equent effectuer l'int\'egration en prenant pour limite le rectangle au lieu du triangle, et divisant par deux; mais on obtient ainsi la quantit\'e 
$$J = \frac12 \int_a^b dx \int_a^b [\varphi(x)\psi(y) - \varphi(y) \psi(x) ]^2 dy,$$
qui se ram\`ene, d'elle m\^eme \`a des int\'egrales simples. 

Les consid\'erations que vous exposez dans votre lettre du 23 Octobre m'ont beaucoup int\'eress\'e et me semblent excellentes; afin de remplir vos intentions je demanderai \`a mon confr\`ere M\mr Darboux d'en prendre connaissance et s'il partage mon sentiment comme j'ai tout lieu de le penser, je le prierai de publier votre lettre soit dans les Annales de l'\'Ecole Normale soit dans le Bulletin des Sciences Math\'ematiques.

Veuillez agr\'eer Monsieur, la nouvelle assurance de mes sentiments de bien sinc\`ere sympathie et de haute estime 

Ch.\ Hermite 

Paris 16 Novembre 1888
\end{quote}

\subsection{November 22, 1888}

This is letter 773 in the archive. This letter refers to a previous letter sent by Gomes Teixeira to Hermite, with a scientific article. It establishes an inequality which can be used to provide an upper bound for integrals of the type 
$$\int_a^b \varphi(x) \psi(x) dx.$$
Teixeira compares this result with an inequality by Chebyschev, which provides a lower bound for integrals of the same type. Hermite makes a small correction in a coefficient and offers publication. The paper mentioned in this letter is \cite{GT52}.

\begin{quote}

Monsieur,

Vous \^etes parvenu \`a des cons\'equences qui m'on paru tr\`es int\'e\-ressan\-tes en rapprochant du Th\'eor\`eme de M.\ Tchebichev, votre relation
$$\int_a^b \varphi^2(x) dx . \int_a^b \psi^2(x) dx > \left[ \int_a^b \varphi(x) \psi(x) dx\right]^2.$$
Je m'en suis entretenu avec M\mr Darboux, en lui donnant communication de votre derni\`ere lettre ou vous les exposez avec plusieurs applications, et je pense point qu'il vous sera d\'esagr\'eable que cette lettre soit publi\'e dans le Bulletin des Sciences Math\'e\-ma\-tiques, comme celle que vous m'aviez d\'ej\`a fait l'honneur de m'adresser. Les \'epreuves vous seront envoy\'es \`a Porto pour que vous puissiez les revoir et les corriger; vous me permettrez de vous proposer de remplacer votre int\'egrale double 
$$\int_a^b dx \int_x^b [\varphi(x) \psi(y) - \varphi(y)\psi(x) ]^2 dy$$
par celle-ci
$$\frac12 \int_a^b dx \int_x^b [\varphi(x) \psi(y) - \varphi(y)\psi(x) ]^2 dy$$
afin qu'on reconnaisse imm\'ediatement que vous partez, comme vous le dites, d'une identit\'e. 

En vous renouvelant, Monsieur, l'assurance de ma haute estime et de mes sentiments d\'evou\'es, 

Ch.\ Hermite

Paris, 22 Novembre 1888

\end{quote}

\subsection{March [17?], 1890}

This is letter 771 in the archive, and it refers to the paper \cite{GT} by Gomes Teixeira, which provides a new proof of a result in Hermite's {\it Cours d'Analyse}, further developed in one of Hermite's articles in the JSMA, namely \cite{HJ2}. Hermite compliments this proof, stating he will actually use it in one of his lessons, and gives yet another proof of the result.

\begin{quote}

Monsieur, 

J'ai lu avec le plus grand plaisir dans le dernier n.$^{\rm o}$ du Journal des Sciences Math\'ematiques et Astronomiques la belle et savante m\'ethode que vous avez expos\'e pour obtenir l'int\'egrale 
$$\int_0^\pi \cot (x-a-ib)dx$$
en la rattachant a l'expression g\'en\'erale 
$$\int_0^\pi \frac{f'(x)dx}{1+f^2(x)}$$ 
et a la notion de l'indice de Cauchy.

Vous avez ainsi obtenu une application tr\`es int\'eressante de la th\'eorie du grand g\'eom\`etre, et je me propose de la mettre \`a profit dans une de mes le\c cons. A cette occasion, et pour vos \'etudes, permettez moi de vous indiquer un autre proc\'ed\'e qui conduit \`a la valeur de cette int\'egrale au moyen de la relation trigonom\'etrique 
$$\cot nx = \sum \frac1n \cot\left( x+\frac{k\pi}n\right), (k=0, 1, 2, \ldots, n-1)$$ 

A cet effet, je change d'abord $x$ en $x-a-ib$, puis je fais cro\^itre ind\'efiniment l'entier $n$. 
En posant $\frac\pi n = dx$, vous voyez que le deuxi\`eme membre a ainsi pour limite l'int\'egrale d\'efinie qu'il s'agit d'obtenir:
$$\frac1\pi \int_0^\pi \cot(x-a-ib)dx,$$
de sorte que l'on a pour $n$ infiniment grand:
$$\int_0^\pi \cot (x-a-ib) dx = \pi \cot n(x-a-ib)$$
Cela \'etant, la formule:
$$\cot n(x-a-ib) = \frac1{i} \frac{e^{2in(x-a-ib)}+1}{e^{2in(x-a-ib)}+1}$$
donne imm\'ediatement dans cette hypoth\`ese $+i$ ou $-i$ suivant que $b$ est positif ou n\'egatif et par cons\'equent le r\'esultat cherch\'e. 

Peut-\^etre pensez vous comme moi Monsieur qu'il est bon dans l'enseignement de donner sous plusieurs points de vue la d\'emons\-tration des propositions importantes, c'est le but que j'ai en vue. Je remarquerai encore que la relation dont j'ai fait usage:
$$\cot nx = \sum \frac1n \cot \left( x+\frac{k\pi}n \right),$$
r\'esulte imm\'ediatement de la d\'ecomposition en \'el\'ements simples de la fonction $\cot nx$, les r\'esidus pour $x=-\frac{k\pi}n$ , c'est \`a dire les valeurs que prend dans cette supposition le produit 
$$\left(x+\frac{k\pi}n\right) \cot nx = \frac{\left( x+ \frac{k\pi}n \right) \cos nx}{\sin nx}$$
\'etant tous \'egaux a $\frac1n$.

Une derni\`ere remarque; l'int\'egrale $J =  \int_0^\pi \cot(t-z) dt$, consid\'er\'ee comme fonction de la variable $z$ offre une coupure repr\'esent\'ee par l'\'equation, $t-z = k\pi$, ou $t$ varie de z\'ero \`a $\pi$, $k$ \'etant un entier arbitraire, cette coupure est dans l'axe des abcisses. Dans le demi plan au dessous de cet axe, l'int\'egrale est une fonction continue uniforme de $z$ et la relation $D_z J = 0$ qu'on obtient imm\'ediatement montre qu'elle est constante. Pour l'obtenir je suppose $z=x+iy$, et en attribuant \`a $y$ successivement des valeurs infiniment grandes, positives et n\'egatives, je trouve encore:
$$J=\int_0^\pi idt, \qquad J=-\int_0^\pi idt,$$
c'est \`a dire, $J=i\pi$ et $J=-i\pi$. 

Veuillez agr\'eer Monsieur, la nouvelle assurance de mes sentiments bien d\'evou\'es

Ch.\ Hermite 

J'abr\`ege un peu le calcul que vous avez donn\'e en partant de l'identit\'e 
$$\cot (x-g) -\cot(x+g) = \frac{ 2\sin g \cos g }{\sin^2 x-\sin^2 g},$$
d'ou l'on tire:
$$\int [\cot (x-g) -\cot(x+g)] dx = \int \frac{ 2\sin g \cos g\, dx}{\sin^2 x-\sin^2 g}$$
Puis si l'on pose $z=\tan g \tan x$
$$\int [\cot (x-g) -\cot(x+g)] dx = \int \frac{dz}{1-z^2}$$
Il suffit ensuite de changer $x$ en $x-a$ et de faire $g=ib$. 

Paris [17?] Mars 1890
\end{quote}

\subsection{July 9, 1890}

This is letter 774 in the archive. This letter refers to article \cite{GT6}, in which a complex integral is studied in order to provide results about series of trigonometric functions. There is also a reference to a text sent to Hermite by Jos\'e Pedro Teixeira (1857--1925), a professor at the Polytechnic Academy of Porto.

\begin{quote}

Paris 9 Juillet 1890

Monsieur, 

Le travail dont vous avez bien voulu me donner communication m'a sembl\'e excellent, et j'esp\`ere que je n'aurai pas \'et\'e contre votre vos intentions en le faisant lire \`a M\mr Darboux qui m'a exprim\'e le d\'esir de le publier dans le Bulletin des Sciences Math\'ematiques. J'ose pensez que vous donnerez votre consentement \`a cette publication qui est un t\'emoignage de mes sentiments de sympathie \`a votre \'egard et de notre estime pour votre talent. Vous voudrez bien aussi Monsieur m'excuser de ne vous soumettre aucune remarque sur votre analyse, les devoirs universitaires sont tellement on\'ereux \`a ce moment, \`a cause du grand nombre d'examens dont je suis oblig\'e \`a la Sorbonne, que le courage et le temps me manquent pour pouvoir m'occuper de sujets \'elev\'es du calcul int\'egral. 

En saisissant cette occasion pour remercier M\mr Jos\'e Pedro Teixiera du don de son m\'emoire sur les fonctions elliptiques, et en vous exprimant toute ma satisfaction de voir cette th\'eorie cultiv\'ee avec autant de succ\`es, je vous prie Monsieur de recevoir la nouvelle assurance e ma plus haute estime et de mes sentiments bien d\'evou\'es.

Ch.\ Hermite
\end{quote} 

\subsection{November 9, 1891}

This is letter 775 in the archive. The letter refers to article \cite{GT7}, which is a study of $\mathcal{p}(u)$, Weierstrass's elliptic function. 

\begin{quote}

Paris 9 Novembre 1891

Monsieur, 

C'est avec le plus grand plaisir que j'ai lu l'analyse si simple et si \'el\'egante dont vous m'avez donn\'e communication.  Votre m\'ethode me semble avoir une des voies les meilleures pour entrer dans la th\'eorie des fonctions elliptiques, et dans le d\'esir que l'on profite du fruit de vos efforts, j'ai pens\'e que vous voudriez bien consentir \`a la publication de votre lettre dans le Bulletin des Sciences Math\'ematiques. Il me fait un agr\'eable devoir de vous exprimer le sentiment de mon confr\`ere M\mr Darboux qui appr\'ecie comme moi votre beau talent, et je saisis Monsieur cette occasion pour vous renouveler l'expression de ma haute estime et de ma bien sinc\`ere sympathie

Ch.\ Hermite.

\end{quote} 

\subsection{May 24, 1892}

This is letter 776 in the archive. In it, Hermite thanks Gomes Teixeira for having sent his {\em Curso de an\'alise}, and makes very specific remarks about it, complementing the clarity and elegance of the text, which  denotes a careful reading. Hermite also promises to send a note for publication, which is sent with letter 777. 

\begin{quote}

Paris 24 Mai 1892

Monsieur,

Les trois volumes de votre Cours d'Analyse Infinit\'esimale que vous m'avez envoy\'es pour \^etre offerts \`a l'Acad\'emie des Sciences, lui ont \'et\'e pr\'esent\'es avec \'eloges par M\mr le Secr\'etaire Perp\'etuel dans la s\'eance du lundi dernier et vous verrez la mention de cette pr\'esentation dans le Compte-rendu de cette s\'eance. Vous ne pouvez douter Monsieur du grand int\'er\^et avec le quel j'ai lu la seconde partie de votre Calcul int\'egral, et tout particuli\`erement la th\'eorie des fonctions doublement p\'eriodiques et des fonctions elliptiques. Il me parait difficile d'atteindre dans l'exposition un plus grand degr\'e de simplicit\'e, de clart\'e et d'\'el\'egance, et je ne puis douter que tous les amis de la science, les g\'eom\`etres comme les commen\c cants ne vous rendent le m\^eme t\'emoignage. J'aurais bien d\'esir\'e vous donner autrement qu'en paroles une preuve de mes sentiments de bien grande estime, en vous envoyant un article pour le journal que vous dirigez avec tant de succ\`es et j'avais song\'e dans ce but \`a r\'ediger une recherche sur les nombres de Bernoulli qui a fait le sujet d'une de mes le\c cons. Mais je suis dans la n\'ecessit\'e d'ajourner ce travail, \'etant \`a cette \'epoque de l'ann\'ee surcharg\'e d'[ouvrage?]; il me faut attendre un moment ou j'aurai plus de libert\'e, et vous voudrez bien me faire cr\'edit d'un peu de temps. 

En vous demandant de prendre note de mon engagement, je saisis Monsieur cette occasion pour vous renouveler l'assurance de ma plus haute estime et celle de mes sentiments bien sinc\`erement d\'evou\'es 

Ch.\ Hermite
\end{quote}

\subsection{December 2 or 4, 1892}

This is letter 777. Even though it contains a reference to a note Hermite sends for publication, the archived material does not include it. As it happens, there is a small article by Hermite in the JSMA, dated about one year and a half after this letter, namely, \cite{HJ4} (it can also be found in \cite[pp.\ 512--513]{HOC}). As there are no more publications by Hermite in the JSMA, other than the three identified above and this one, we believe this note to be the one mentioned in this letter. 

The paper concerns a property of elliptic functions and elliptic integrals, providing a proof for a formula about addition of arguments in Weierstrass's elliptic function $\mathcal{p}(x)$. 

\begin{quote}
Paris 4 D\'ecembre 1892

Monsieur,

J'ai un double remerciement \`a vous adresser; je vous rend gr\^aces d'abord pour votre communication dont j'ai donn\'e connaissance, en pensant remplir votre intention, \`a mon confr\`ere de l'Acad\'emie des Sciences, M\mr Darboux, qui l'a accueillie avec empressement pour la publier dans le Bulletin des Sciences Math\'ematiques. Il me faut ensuite vous t\'emoigner ma vive gratitude d'avoir bien voulu joindre votre nom estim\'e de tous les g\'eom\`etres \`a ceux de mes amis math\'ematiques qui ont d\'esir\'e m'offrir un t\'emoignage de sympathie \`a l'occasion de mes 70 ans. La sympathie est r\'eciproque, et vous y avez un droit particulier, en raison du talent et du z\`ele que vous consacrez a votre Journal des Sciences Math\'ematiques et Astronomiques, et le service \'eminent que vous rendez ainsi \`a la science de votre pays. Peut-\^etre n'avez vous pas oubli\'e que je vous avais promis une note que devrait concerner le th\'eor\`eme de Staudt sur les nombres de Bernoulli, mais mon travail a pris une autre direction et c'est avec un autre sujet que je m'acquitterai envers vous. Permettez moi Monsieur de vous la transcrire ci-apr\`es et de saisir cette occasion pour vous renouveler l'assurance de ma haute estime et de mes sentiments bien d\'evou\'es 

Ch.\ Hermite 

Paris 2 D\'ecembre 1892

\end{quote} 

\subsection{January 28, 1896} 

This is letter 30 in the archive. It contains a request from Hermite to add Gomes Teixeira's name to a letter in support of {\it Acta Mathematica}, a journal founded in 1882 by G\"osta Mittag-Leffler, which was going through difficult financial problems at the time. The request is so urgent that an absence of response is considered an affirmation of support. It is clearly a text that was sent to many scientists at the time: a very similar letter can be found in \cite{Co}, with the same date, addressed to Markov. This paper also transcribes the text of the circular letter from paper \cite{No}, where the names of 30 signatories, among which Gomes Teixeira, can be found. It adds that there are 370 more signatures, in a total of 400 mathematicians. 

\begin{quote} 

Monsieur, 

Permettez moi de faire un appel \`a votre bonne obligeance en sollicitant votre concours en faveur d'un int\'er\^et scientifique au quel j'esp\`ere vous ne refuserez pas votre sympathie. 

Les Acta Mathematica traversent en ce moment une crise qui met leur existence en grand p\'eril; la commission du budget de la di\`ete su\' edoise a manifest\'e l'intention de supprimer compl\`etement l'allocation, apr\`es l'avoir diminu\'ee de moiti\'e l'ann\'ee derni\`ere, qu'elle avait accord\'e au Journal, depuis sa fondation. Ce n'est pas seulement en France que les \'economies sont [?] et [r\'eclament?] sans piti\'e, et que les assembl\'ees politiques se soucient peu des g\' eom\`etres. M$^{\text r}$ Mittag-Leffler et les Acta ne sont point cependant sans d\'efenseurs, par une circonstance heureuse l'un des plus [d\'edi\'es?] est le pr\' esident m\^eme de la commission du budget \`a la di\`ete. Ce haut patronage a \'emis l'avis qu'une manifestation venue de l'\'etranger, apportant au nom de savants autoris\'es un t\'emoignage public de l'importance des Acta et du m\'erite \'elev\'e de son fondateur, d\'eciderai certainement un vote favorable de l'assembl\'ee. L'Acad\'emie des Sciences de Paris est entr\'ee dans ces vues avec empressement et \`a l'unanimit\'e; elle a pris l'initiative d'une adresse dont le 90$^{\text{me}}$ anniversaire de naissance de M$^{\text r}$ Mittag-Leffler offrait l'occasion, pour lui exprimer ses sentiments de sympathie, en lui faisant don de son portrait. En m\^eme temps, elle se propose d'adresser une circulaire aux g\'eom\`etres pour obtenir leur adh\'esion, et j'ai re\c{c}u Monsieur la mission de vous soumettre la circulaire et l'adresse, en vous demandant de permettre que votre nom soit joint \`a ceux de mes confr\`eres de l'Acad\'emie des Sciences et des math\'ematiciens qui nous ont d\'ej\`a donn\'e leur concours. Comme le temps nous est strictement mesur\'e, je prends la libert\'e de vous demander une r\'eponse, pour le seul cas o\`u elle ne r\'epondrait pas a notre atteinte, et par une simple carte postale, dans le moindre d\'elai que vous sera possible, en consid\'erant que votre acceptation nous est acquise si je ne re\c{c}ois pas d'avis contraire. 

Dans l'esp\'erance que je n'aurai pas fait appel en vain \`a votre concours, dans cette circonstance, je saisis Monsieur cette occasion pour renouveler l'assurance de ma plus haute estime et de mes sentiments bien d\'evou\'es.

Ch.\ Hermite 

Paris 28 Janvier 1896

\end{quote} 

\section{Conclusions} 

As we consider these texts, we would like to draw attention to two aspects which become apparent in the relationship between Hermite and Gomes Teixeira. 

This correspondence shows a notable personality trait of Hermite: his desire to continuously support and encourage a younger colleague, in a country that did not have a very strong international presence in mathematics. Even if we take into account the courteous formalities in the text, one can nevertheless say that there is some genuine interest and encouragement from Hermite's part towards Gomes Teixeira. 

It is likewise remarkable that Gomes Teixeira started this correspondence at a very early age: he was still 21 years old when he first contacted Hermite, one of the leading names in mathematics, asking him to participate in a journal that was, at the time, only a project. It is equally remarkable that Hermite did send a paper in reply. 

The texts of the letters presented here may not show the same depth as the ones between Hermite and other of his numerous correspondents. To take only two examples, the correspondence with Mittag-Leffler (which can be found at the {\em Cahiers du s\'eminaire d'histoire des math\'ematiques}, see \cite {HM1}, \cite {HM2} and \cite{HM3}) denotes a more personal relationship, and the letters to Jacobi (some of which can be found at the {\em Journal f\"ur die reine und angewandte Mathematik}, see \cite{HJ4}) show greater mathematical depth. However, this correspondence was maintained throughout 21 years, from 1872 to 1892, and tells of a persistent relationship maintained between the two mathematicians. Moreover, the nine letters addressed to Gomes Teixeira about the homage to Hermite in 1892 (not included in this paper) show that this proximity was acknowledged by the mathematical community of their time. 

This analysis of the correspondence from Hermite to Gomes Teixeira is, to the best of our knowledge, the first in-depth study of the material in the archive of the correspondence received by Gomes Teixeira. It focuses only on 18 letters, out of the more that 2000 that are extant in this archive. We hope to continue this work in following publications, focusing on other aspects of Gomes Teixeira's correspondence. 

\section{Acknowledgements} 

The author wishes to thank Catherine Goldstein for the revision of the transcriptions, as well as the other members of this project, namely Henrique Leit\~ao, Jo\~ao Queir\'o and Ant\'onio Leal Duarte for many contributions and exchanges of opinions. 

Thanks also to the Archive of the University of Coimbra for the kind cooperation, in the access to the letters, and to the Funda\c c\~ao para a Ci\^encia e Tecnologia, for funding, through Project FCT UID/HIS/00286/2013.

\end{document}